\documentclass{article}
\usepackage{amssymb}

\usepackage{graphicx}
\usepackage{amsmath}
\newenvironment{proof}[1][Proof]{\textbf{#1.} }{\hfill  $\Box$}
\newtheorem{theorem}{Theorem}

\newtheorem{corollary}[theorem]{Corollary}

\newtheorem{definition}[theorem]{Definition}

\newtheorem{lemma}[theorem]{Lemma}

\newtheorem{proposition}[theorem]{Proposition}

\begin{document}

\title{\textbf{Generalized Teichm\"{u}ller space of non-compact 3}$-$\textbf{%
manifolds and Mostow rigidity}}
\author{Charalampos Charitos and Ioannis Papadoperakis}
\maketitle

\begin{abstract}
Consider a 3$-$dimensional manifold $N$ obtained by gluing a finite number of ideal
hyperbolic tetrahedra via isometries along their faces. By varying the
isometry type of each tetrahedron but keeping fixed the gluing pattern 
 we define a space $\mathcal{T}$ of complete hyperbolic metrics on $N$
with cone singularities along the edges of the tetrahedra.
We prove that $\mathcal{T}$ is homeomorphic to a Euclidean space and we
compute its dimension. By means of examples, we examine if the elements of $%
\mathcal{T}$ are uniquely determined by the angles around the edges of $N.$

\end{abstract}

\section{Introduction}

In \cite{[ChAth]}, \cite{[ChAth1]} spaces $X$ which are called \textit{ideal
simplicial complexes, }are considered. These spaces $X$ are obtained by
gluing along their edges finitely many ideal hyperbolic triangles. The
Teichm\"{u}ller space $\mathcal{T}(X)$ of $X$ is defined and parametrized
via the shifts parameters.

In the present work we consider orientable, compact 3$-$manifolds with
non-empty boundary $\partial M.$ The interior $Int(M)$ of $M$ always has a
triangulation $\mathcal{D}$ by ideal tetrahedra. Fixing $\mathcal{D}$ we
define \textit{ideal hyperbolic structures} with axial singularities on $M,$
as well as, the \textit{generalized Teichm\"{u}ller space} $\mathcal{T}_{\mathcal{D}}(M)$
of $M.$ The 2$-$skeleton of $\mathcal{D},$ equipped with a hyperbolic
structure induced from an element of $\mathcal{T}_{\mathcal{D}}(M),$ is a 2$-$dimensional
ideal simplicial complex $X$ and we show that the shift parameters\textit{\ }
which parametrize $\mathcal{T}(X)$ also parametrize $\mathcal{T}_{\mathcal{D}}(M).$ Thus $%
\mathcal{T}_{\mathcal{D}}(M)$ is homeomorphic to a Euclidean space $\mathbb{R}^{d}$ and
we prove that $d$ is equal to the number of edges minus the number of
vertices of $\mathcal{D}.$

By Mostow rigidity theorem, if $h,$ $h^{\prime }$ are ideal hyperbolic
structures on $M$ and the angle around each edge of $\mathcal{D}$ is equal
to $2\pi ,$ then $h,$ $h^{\prime }$ represent the same element in $\mathcal{T%
}(M).$ In this work we give examples of $3-$manifolds equipped with ideal
hyperbolic structures and we show that all these structures are parametrized
by the angles around the edges of $\mathcal{D}.$ An interesting problem for
further investigation, is to consider ideal hyperbolic structures on $M,$
i.e. complete metrics $h$ in the interior of $M$ of curvature $\leq $ $-1,$
and examine if the angles around the edges of $\mathcal{D}$ uniquely
determine $h$ as an element of $\mathcal{T}_{\mathcal{D}}(M).$

\section{Definitions and Preliminaries}

In his pioneering work \cite{[ThurstonNotes]}, Thurston constructed a
hyperbolic structure on the complement of certain knots by realizing them as
a union of finitely many ideal hyperbolic tetrahedra. In the present paper,
inspired from Thurston's work and from the work of other mathematicians, see
for example \cite{[fmp]}, \cite{[FP]}, we glue a finite number of ideal
hyperbolic tetrahedra and we consider, in the resulting manifold, hyperbolic
structures in a broader sense.

\begin{definition}
\label{origin}Assume that $M$ is a compact, orientable 3$-$manifold with $%
\partial M\neq \emptyset .$ A\textit{\ topological ideal triangulation} of $%
M $ consists of two finite sets $\mathcal{D}$ and $\mathcal{F}$ which
satisfy the following two conditions:

(1) Each element $\Delta \in \mathcal{D}$ is a standard tetrahedron and each
element $f\in \mathcal{F}$ is a simplicial homeomorphism $f:A\rightarrow B,$ where $A$
and $B$ are triangular faces of two tetrahedra $\Delta $ and $%
\Delta ^{\prime }$ of $\mathcal{D}.$ The elements of $\mathcal{F}$ are
called \textit{gluing maps }and they are orientation reversing simplicial
maps. Furthermore, for each face $A$ of a tetrahedron $\Delta \in \mathcal{D}%
,$ there exists precisely one $f\in \mathcal{F}$ and a face $B$ of some tetrahedron $%
\Delta ^{\prime }\in \mathcal{D},$ such that $f$ maps $A$ onto $B$ or $B$ onto $A.$
 
(2) If $Y$ is the quotient space of the disjoint union of all tetrahedra in $%
\mathcal{D}$ by the relation which identifies any two points $x\in A$ and $%
y\in B$ by a map $f\in $ $\mathcal{F}$ satisfying $f(x)=y$ and if we remove
from $Y$ all vertices of tetrahedra then we obtain a space homeomorphic to
the interior $Int(M)$ of $M.$

The subdivision of $Int(M)$ into tetrahedra of $\mathcal{D}$ with its
vertices deleted, will be called \textit{\ topological ideal triangulation }%
of $M$ and will be also denoted by $\mathcal{D}.$ Each tetrahedron $\Delta
\in $ $\mathcal{D}$ will be called an \textit{\ ideal tetrahedron}. A face
(resp. an edge) of some $\Delta $ in $\mathcal{D}$ will be called a face
(resp. an edge) of $\mathcal{D}.$ The deleted vertices of $Y$ will be called
ideal vertices of $Int(M)$ or vertices of $\mathcal{D}.$
\end{definition}

\noindent \textbf{Remark} A more accurate picture is obtained, rather than
by removing the vertices, by truncating the tetrahedra; that is, by removing
an open neighborhood of each vertex of tetrahedra. Then we recover, not only 
$Int(M),$ but the whole $M$ by these truncated tetrahedra.\smallskip

Henceforward, for each manifold $M$ we will denote by $M^{o}$ its interior.
We shall be interested in metrics $h$ on $M^{o}$ which are obtained in the
following manner: each tetrahedron $\Delta \in $ $\mathcal{D}$ is equipped
with a metric which makes it isometric to an ideal hyperbolic tetrahedron.
These tetrahedra are glued among them along isometries 
and thus $M^{o}$ is equipped naturally with the length metric. The
subdivision of $M^{o}$ into ideal hyperbolic tetrahedra will be called \textit{%
\ hyperbolic ideal triangulation }of $M$ and will also be denoted by $%
\mathcal{D}.$ The length metric $h$ on $M^{o}$ will be called an \textit{%
ideal metric}.

If $e$ is an edge of $\mathcal{D},$ we denote by $\theta _{h}(e)$ the sum of
all dihedral angles formed by the faces of ideal hyperbolic tetrahedra which
have $e$ as a common edge. Then we distinguish two cases:

$(1)$ If $\theta _{h}(e)\neq 2\pi ,$ the edge $e$ is called \textit{singular}
or \textit{axis}.

$(2)$ If $\theta _{h}(e)=2\pi $ the edge $e$ is called \textit{regular}.

\noindent$\theta _{h}(e)$ will be called the \textit{angle around the edge }$%
e.$

We shall henceforth assume that the metric $h$ on $M^{o}$ is complete. This
metric $h$ has singularities along the axes of $M^{o}$ which will be called 
\textit{axial singularities. }Besides the axes of $M^{o},$ the curvature of $%
h$ is constant, equal to $-1.$ Such a complete metric $h$ will be called an 
\textit{ideal structure} on $M,$ with respect to $\mathcal{D}.$ In what
follows the topological ideal triangulation $\mathcal{D}$ will be fixed so
the specification ``with respect to $\mathcal{D}"$ is omitted.

The completeness of $h$ imposes some restrictions on the gluing maps, which
can easily be described in terms of a geometrical property at the \textit{%
ideal vertices or cusps} of $M^{o}.$ (When an ideal metric $h$ is considered
on $M^{o},$ an ideal vertex of $M^{o}$ will be also referred as \textit{cusp}
of $M^{o}).$ For each cusp $v$ of $M^{o},$ we can associate a natural
foliation of a subset of $M^{o}$ (a ``neighborhood'' of $v$). The definition
is as follows. Consider an ideal hyperbolic tetrahedron $\Delta $ in $%
\mathcal{D},$ having $v$ as one of its ideal vertices. Consider a foliation
of a horoball neighborhood of $v$ in $\Delta ,$ whose leaves are pieces of
horodiscs which are centered at $v.$ Then, $(M^{o},h)$ is complete as a
metric space if and only if, the horodiscs on each ideal hyperbolic
tetrahedron abutting at $v,$ fit together properly so that they form a
product foliation $K\times \{t\},$ $t\in \lbrack 0,\infty )$ defined in a
``horoball neighborhood'' $V$ of $v$ in $M^{o}.$ The complete metric $h$
defined on $M^{o},$ induces on every fiber $K_{t}=K\times \{t\},$ which is a
closed surface, a Euclidean structure $h_{t}$ with conical singularities,
see \cite{[Troyanov0]} for the precise definition and for a thorough
discussion of structures $h_{t}.$ In fact, each $K_{t}$ is naturally
triangulated by the horospherical section of the ideal tetrahedra. The
conical singularities arise exactly at the points where the singular edges
(axes) intersect $K_{t}.$ These sections are Euclidean triangles and, since $%
h$ is complete, they are glued by isometries. Obviously, for each $t\in
\lbrack 0,+\infty ),$ $h_{t}$ is a rescaling of the metric $h_{0}.$ A
surface $S$ which coincides with some $K_{t},$ $t\in \lbrack 0,+\infty )$
will be referred to as the \textit{geometrical link} of $v,$ with respect to 
$h.$

\section{The Teichm\"{u}ller space of the 2$-$skeleton of a hyperbolic ideal
triangulation}

Consider the topological ideal triangulation $\mathcal{D}$ of $M$ and let $%
\mathcal{D}^{(2)}$ be the 2$-$skeleton of $\mathcal{D}.$ Let $h$ be an ideal
structure on $M^{o}.$ With respect to $h,$ every tetrahedron of $\mathcal{D}$
with its vertices deleted becomes an ideal hyperbolic tetrahedron and so
every face of $\mathcal{D}^{(2)}$ is isometric to an ideal hyperbolic
triangle. Denote by $X=|\mathcal{D}^{(2)}|$ the support of $\mathcal{D}%
^{(2)} $ and let $\overline{h}$ be the metric induced by $h$ on $X.$ Then $X$
equipped with $\overline{h},$ is an \textit{ideal 2}$-$\textit{dimensional
simplicial complex} in the sense of \cite[Def. 3.1]{[ChAth1]}. An \textit{%
edge }of $X$ is a 1$-$simplex of $\mathcal{D}^{(2)}$ and it is isometric to
a line. A \textit{face} of $X$ is a 2$-$simplex of $\mathcal{D}^{(2)}$ and
it is isometric to a hyperbolic ideal triangle. The deleted vertices of $%
\mathcal{D}^{(2)}$ are called \textit{cusps} of $X.$

Let $v$ be a cusp of $X.$ The link $\Gamma =\Gamma (v)$ of $v$ in $X$ is a
simplicial graph embedded in $X,$ which is defined by taking one vertex on
each half-edge of $X$ abutting on $v$ and then joining two such vertices by
an edge contained in a face of $X.$ Notice also that there exists a closed
neighborhood $V$ in $X$ ``centered'' at $v$ which has a natural structure of
a geometric cone $v\cdot \Gamma -\{v\}$ (see Def. 2.2 in \cite{[ChAth]}). We
will be calling $V$ a \textit{neighborhood} of $v.$ We have the following
lemma.

\begin{lemma}
\label{induced metric} The ideal metric $h$ on $M^{o}$ is complete if and
only if the induced metric $\overline{h}$ is complete on $X.$
\end{lemma}

\begin{proof}
First assume that $h$ is complete. Since $X$ is a closed subset of $M^{o}$
we have that the induced metric $\overline{h}$ on $X$ is complete. Assume
now that $\overline{h}$ is complete. Then, for each cusp $v$ of $X$ the
horocycles on each ideal triangle which have $v$ as an ideal vertex, fit
together properly so that they form, in a neighborhood of $v,$ a connected
graph whose edges are horocycle segments, see Proposition 3.4.18 in \cite
{[Thurstonbook]}. Actually, this proposition is proven for cusped surfaces
but the same method of proof applies for ideal 2$-$dimensional simplicial
complexes. This implies that the horospherical sections of every ideal
hyperbolic tetrahedron on a neighborhood of $v,$ fit together forming a
closed surface which is the geometrical link of $v.$ From the discussion in
the section above, we deduce that $h$ is complete.
\end{proof}

By considering various ideal structures $h$ on $M$ we obtain various \textit{%
ideal hyperbolic structures} $\overline{h}$ on $X$ i.e. complete metrics
such that $(X,\overline{h})$ is a local $CAT(-1)$ space, see Prop. 1.4 of 
\cite{[ChAth]}. This leads us to consider the Teichm\"{u}ller space $%
\mathcal{T}(X)$ of $X$ and relate it to the generalized Teichm\"{u}ller
space $\mathcal{T}_{\mathcal{D}}(M)$ of $M$ which will be defined in the next section.

We recall the definition of $\mathcal{T}(X)$ (see Def. 2.1 in \cite{[ChAth]}%
).

\begin{definition}
\label{defteich}The Teichm\"{u}ller space $\mathcal{T}(X)$ of $X$ is the set
of equivalence classes of ideal hyperbolic structures on $X;$ such two
structures $h$ and $h^{\prime }$ are considered equivalent if there is a
homeomorphism $F:X\rightarrow X$ which preserves each edge and each face of $%
X$ and which satisfies $F^{\ast }(h)=h^{\prime },$ where $F^{\ast }(h)$
denotes the pull-back of the metric $h$ via $F.$ Remark that such an $F$ is
isotopic to the identity map.
\end{definition}

Let $T$ be an ideal hyperbolic triangle. Then $T$ has a distinguished point
which is the \textit{barycentre} of $T.$ Each edge of $T$ is also equipped
with a distinguished point, namely, the foot of the perpendicular drawn from
the barycentre of $T$ to that edge. We shall call this point the \textit{%
centre} of the edge.

There are several ways of describing the topology of $\mathcal{T}(X),$ and
we shall use here the \textit{shift} parameters. Let $X$ be an ideal $2-$%
dimensional simplicial complex equipped with an ideal hyperbolic structure $%
h $ and let $V$ be the set of cusps of $X.$ In order to describe the shift
parameters, we start by choosing once and for all an orientation on each
edge of $X.$ If $T,$ $T^{\prime }$ are two faces of $X$ with $e\subset T\cap
T^{\prime },$ we define the quantity $x_{h}(T,T^{\prime },e)$ as the
algebraic distance on $e$ from the centre $p$ of $e$ associated to $T$ to
the centre $p^{\prime }$ of $e$ associated to $T^{\prime },$ and we call it
the \textit{shift parameter} on the ordered triad $(T,T^{\prime },e).$ More
precisely, if the direction from $p$ to $p^{\prime }$ coincides with the
orientation of $e$ then $x_{h}(T,T^{\prime },e)$ is positive, otherwise it
is negative (see Def. 3.2 in \cite{[ChAth1]}].

Let $\mathcal{B}$ be the set of ordered triads $(T,T^{\prime },e)$ where $T,$
$T^{\prime }$ are triangles of $X.$ The shift parameter defines a map $%
\mathcal{I}:\mathcal{T}(X)\rightarrow \mathbb{R}^{\mathcal{B}},$ by the
formula

\begin{center}
$\mathcal{I}(h)(T,T^{\prime },e)=x_{h}(T,T^{\prime },e).$
\end{center}

The map $\mathcal{I}$ is clearly injective but not necessarily onto. We
equip $\mathcal{T}(X)$ with the topology induced
from the embedding $\mathcal{I}:\mathcal{T}(X)\rightarrow \mathbb{R}^{%
\mathcal{B}}.$ Thus, $\mathcal{T}(X)$ is parametrized by the shift parameters
of the elements of $\mathcal{B}.$ These shift parameters satisfy certain linerar 
equations and therefore  $\mathcal{T}(X)$ is homeomorphic
to a Euclidean space. By induction on the number of triangles
of $X$ we obtain the following Proposition, which is stated without proof in
 \cite{[ChAth]}.

\begin{proposition}
\label{dim polyhedron} Let $d_{0}$ be the number of gluing maps $\phi $
appearing in the construction of $X$ and let $r_{i}$ be the rank of $\pi
_{1}(\Gamma _{i}).$ Then $\mathcal{T}(X)$ is homeomorphic to a Euclidean
space and its dimension is equal to $d_{0}-\sum_{i=1}^{k}r_{i}.\qquad $
\end{proposition}

\section{The generalized Teichm\"{u}ller space of $M$}

In this section we will prove that if $h^{\prime }$ is an ideal structure on 
$M$ such that $\theta _{h^{\prime }}(e)\geq 2\pi $ for each edge $e$ of $%
\mathcal{D}$ then any other ideal structure $h$ on $M$ is hyperbolic in the
sense of Gromov \cite{[Gromov]}. After that, we will define the generalized
Teichm\"{u}ller space of $M$ and we will compute its dimension.

A map $f:X\rightarrow Y$ between metric spaces is \textit{Lipschitz} if and
only if there is a constant $K>0$ such that 
\begin{equation*}
d_{Y}(f(x),f(y))\leq Kd_{X}(x,y)\text{ for all }x,y\text{ in }X.
\end{equation*}

\noindent If $f$ is a homeomorphism then $f$ is called \textit{bi-Lipschitz }%
if $f$ and $f^{-1}$ are \textit{Lipschitz.}

Let $h$ be an ideal structure on $M.$ Denote by $(M^{o},h)$ the interior of $%
M$ equipped with the complete metric $h$ and by $(\widetilde{M^{o}},%
\widetilde{h})$ the universal covering of $M^{o}$ equipped with a metric $%
\widetilde{h}$ such that the covering projection $\pi :(\widetilde{M^{o}},%
\widetilde{h})\rightarrow (M^{o},h)$ is a local isometry. Let $h,$ $%
h^{\prime }$ be two ideal structures on $M.$ We have the following
proposition.

\begin{proposition}
\label{lips2} $(1)$ The identity map $Id:$ $(M^{o},h)\rightarrow
(M^{o},h^{\prime })$ is bi-Lipschitz.

\noindent \noindent $(2)$ There is a lifting $\widetilde{Id}:(\widetilde{%
M^{o}},\widetilde{h})\rightarrow (\widetilde{M^{o}},\widetilde{h^{\prime }})$
which is bi-Lipschitz.
\end{proposition}

\begin{proof}
$(1)$ Let $\Delta $ be a tetrahedron of $\mathcal{D}.$ Obviously $%
Id_{|\Delta }$ sends $(\Delta ,h)$ onto $(\Delta ,h^{\prime }).$
Furthermore, using appropriate coordinates, $Id_{|\Delta }$ preserves the
hyperbolic height on the link of each vertex of $\Delta .$ This immediately
implies that $Id_{|\Delta }$ is bi-Lipschitz.

The spaces $(M^{o},h)$ and $(M^{o},h^{\prime })$ are geodesic. Let $x,y$ be
two arbitrary points of $M^{o}$ and let $\gamma _{h^{\prime }}[x,y]$ be a
geodesic segment joining $x$ and $y,$ which realizes the distance of $x,y$
with respect to $h^{\prime }.$ Let $x_{1}=x,$ $x_{2},...,$ $x_{n}=y$ be a
subdivision of $\gamma _{h^{\prime }}[x,y]$ such that $x_{i},$ $x_{i+1}$
belong to the same ideal tetrahedron of $\mathcal{D}$ for each $i.$ Then $%
h(x,y)\leq \sum_{i}h(x_{i},x_{i+1})$ \noindent and since $Id_{|\Delta }$ is
bi-Lipschitz we have that there exists $K>0$ such that

\begin{equation*}
h(x,y)\leq \sum_{i}h(x_{i},x_{i+1})\leq \sum_{i}Kh^{\prime
}(x_{i},x_{i+1})=Kh^{\prime }(x,y)
\end{equation*}

\noindent \noindent This last inequality proves that $Id$ is Lipschitz.
Similarly we prove that $Id$ is bi-Lipschitz.

\noindent \noindent $(2)$ Obviously $Id$ lifts to Lipschitz homeomorphism $%
\widetilde{Id}:(\widetilde{M^{o}},\widetilde{h})\rightarrow (\widetilde{M^{o}%
},\widetilde{h^{\prime }})$ which is also bi-Lipschitz.\smallskip
\end{proof}

Let now $\mathcal{H}(M)$ be the set of all ideal structures on $M.$

\begin{corollary}
Assume that there exists $h^{\prime }\in \mathcal{H}(M)$ such that that $%
\theta _{h^{\prime }}(e)\geq 2\pi $ for each edge $e$ of $\mathcal{D}.$ Then
for any metric $h\in \mathcal{H}(M)$ the metric space $(\widetilde{M^{o}},%
\widetilde{h})$ is hyperbolic in the sense of Gromov.
\end{corollary}

\begin{proof}
The metric space $(M^{o},h^{\prime })$ is of curvature less than or equal to 
$-1$ i.e. $(M^{o},h^{\prime })$ satisfies locally the $CAT(-1)$ inequality
(see Thm. 3.13 in \cite{[Paulin]}). Therefore $(\widetilde{M^{o}},\widetilde{%
h^{\prime }})$ is a $CAT(-1)$ space (see \cite{[Gromov]}, page 119) which
implies that $(\widetilde{M^{o}},\widetilde{h^{\prime }})$ is hyperbolic in
the sense of Gromov. Now, from Proposition \ref{lips2}, there is a
bi-Lipschitz mapping between $(\widetilde{M^{o}},\widetilde{h})$ and $(%
\widetilde{M^{o}},\widetilde{h^{\prime }}).$ Therefore $(\widetilde{M^{o}},%
\widetilde{h})$ is hyperbolic in the sense of Gromov, see Thm. 2.2 in \cite
{[CDP]}.
\end{proof}

Such an ideal structure $h$ will be referred to as \textit{ideal hyperbolic
structure.} So, the previous lemma asserts that if $M$ admits an ideal
hyperbolic structure then all ideal structures on $M$ are hyperbolic. For
these manifolds $M$ we will define the generalized Teichm\"{u}ller space $%
\mathcal{T}_{\mathcal{D}}(M)$ of $M.$

\begin{definition}
\label{defgemteichspace}The generalized Teichm\"{u}ller space $\mathcal{T}%
(M) $ of $M$ is the set of equivalence classes of ideal hyperbolic
structures on $M;$ such two structures $h$ and $h^{\prime }$ are considered
equivalent if there is a homeomorphism $F:M^{o}\rightarrow M^{o}$ which
preserves each ideal tetrahedron, each face and each edge of $\mathcal{D}$
and which satisfies $F^{\ast }(h)=h^{\prime }.$
\end{definition}

Let $h\in \mathcal{T}_{\mathcal{D}}(M).$ If $\overline{h}$ is the induced metric on $X=|%
\mathcal{D}^{(2)}|$ then, from the definition of spaces $\mathcal{T}_{\mathcal{D}}(M)$ and 
$\mathcal{T}(X)$ and from Lemma \ref{induced metric}, we may immediately
deduce that the mapping $\Psi :\mathcal{T}_{\mathcal{D}}(M)$ $\rightarrow $ $\mathcal{T}%
(X) $ which sends the equivalence class of $h$ in $\mathcal{T}_{\mathcal{D}}(M)$ to the
equivalence class of $\overline{h}$ in $\mathcal{T}(X)$ is a bijection.
Therefore, we may equip $\mathcal{T}_{\mathcal{D}}(M)$ with the topology induced from $%
\Psi ,$ so that $\Psi $ becomes a homeomorphism and $\mathcal{T}_{\mathcal{D}}(M)$
homeomorphic to $\mathbb{R}^{d}$ for some $d.$

We have the following theorem.

\begin{theorem}
\label{dimension}The generalized Teichm\"{u}ller space $\mathcal{T}_{\mathcal{D}}(M)$ of $%
M $ is homeomorphic to $\mathbb{R}^{d},$ where $d$ is equal to the number of
edges minus the number of vertices of $\mathcal{D}.$
\end{theorem}

\begin{proof}
Let $E,$ $V$ and $F$ be the set of edges, vertices and faces, respectively,
of the topological ideal triangulation $\mathcal{D}$ of $M.$

Consider an element $h\in \mathcal{T}_{\mathcal{D}}(M).$ For every cusp $v_{i}$ of $M^{o},$
let $S_{h}^{i}$ be the geometrical link of $v_{i}$ with respect to $h.$ $%
\mathcal{D}$ induces a loose Euclidean triangulation $\mathcal{D}_{e}^{i}$
by Euclidean pseudo-triangles on every $S_{h}^{i},$ i.e. each
pseudo-triangle of $\mathcal{D}_{e}^{i}$ is isometric to a Euclidean
triangle which may have two or three vertices indentified to one point.
Remark also that the triangles of $\mathcal{D}_{e}^{i}$ can be multiply
incident to each other. Denote by $T_{i},$ $A_{i},$ $K_{i}$ the sets of
triangles, edges, vertices of $S_{h}^{i}$ respectively. Denote also by $T,$ $%
A,$ $K$ the sets of all triangles, edges and vertices of $S=\cup
_{i}S_{h}^{i}.$

Now, for the number $d_{0}$ of Proposition \ref{dim polyhedron}, we have
that,

\begin{equation}
d_{0}=3card(F)-card(E)\qquad  \tag{1}  \label{one}
\end{equation}
This follows from the fact that, if an edge $e$ of $X$ belongs to $n$
different faces of $\mathcal{D},$ then the number of gluing isometries that
identify these faces along $e$ is $n-1.$

Now, every edge of $\mathcal{D}$ intersects $S$ in two vertices and every
face of $\mathcal{D}$ intersects $S$ in three edges. Therefore, we
have\smallskip\ that, $card(F)=\frac{card(A)}{3}$ and $card(E)=\frac{card(K)%
}{2}.$ By replacing these relations to \ref{one} we have that,

\begin{equation}
d_{0}=card(A)-\frac{card(K)}{2}\qquad  \tag{2}  \label{two}
\end{equation}
Now, from Proposition \ref{dim polyhedron}, we have that $\dim \mathcal{T}%
(X)=$ $d_{0}-\sum_{i=1}^{k}r_{i}.$ On the other hand, a well known fact from
graph theory asserts that $r_{i}=card(A_{i})-card(K_{i})+1,$ $\forall i.$
Therefore, from relation \ref{two} we have that $\dim \mathcal{T}(X)=\frac{%
card(K)}{2}-card(V)=$ $card(E)-card(V)$ which proves the theorem.
\end{proof}

\noindent \hspace*{-7cm} \includegraphics[scale=0.7]{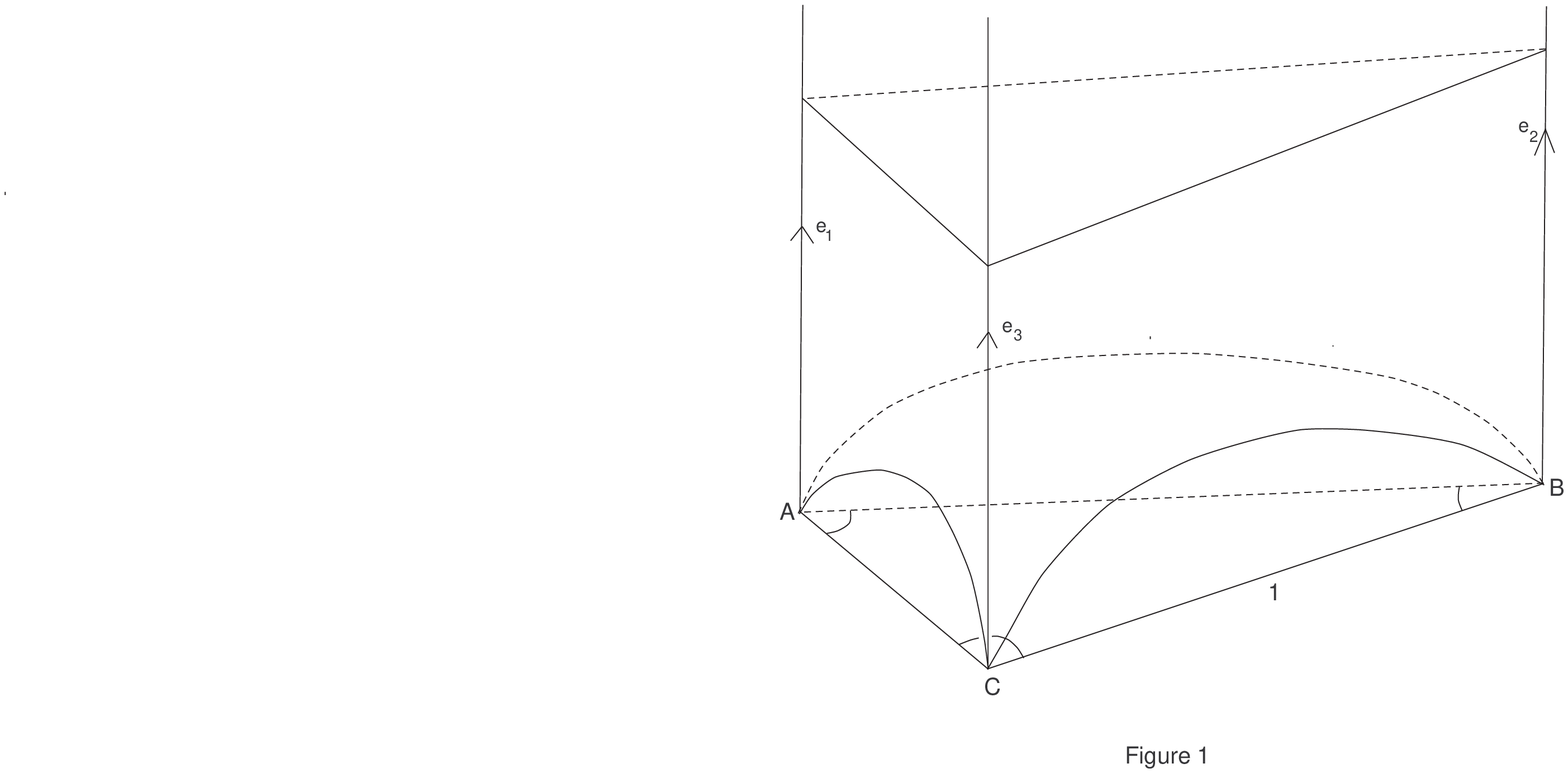}
\section{The angles of axes}

Let $\Delta $ be an ideal hyperbolic tetrahedron in $\mathbb{H}^{3}$ which
has an ideal vertex at $\infty .$ We equip the edges of $\Delta $ with an
orientation such that the edges $e_{1},$ $e_{2},$ $e_{3}$ abutting on $%
\infty $ are oriented towards $\infty .$ Let also $\alpha ,$ $\beta ,$ $%
\gamma $ be the dihedral angles of $\Delta $ corresponding to $e_{1},$ $%
e_{2},$ $e_{3}$ and let  $T_{1}=(C,B,\infty ),$ $T_{2}=(A,C,\infty ),$ 
$T_{3}=(A,B,\infty ),$ see Figure 1.

The tetrahedron $\Delta $ is parametrized by the angles $\alpha ,$ $\beta ,$ 
$\gamma $ which satisfy the relations $0<\alpha ,$ $\beta ,$ $\gamma <\pi $
and $\alpha +$ $\beta +$ $\gamma =\pi .$ Therefore, each $(\alpha ,\beta
,\gamma )$ determines a unique point in the interior of a triangle $T\subset
(0,\pi )^{3}$ whose vertices are the points $(\pi ,0,0),$ $(0,\pi ,0),$ $%
(0,0,\pi ).$ On the other hand, the boundary $\partial \Delta $ of $\Delta $
equipped with the hyperbolic metric from $\Delta ,$ say $h,$ is an ideal 2$-$%
dimensional simplicial complex which is homeomorphic to the sphere $\mathbb{S%
}^{2}-\{3$ points\}. If we consider the shift parameters $\varkappa
_{1}=x_{h}(T_{3},T_{1},e_{1}),$ $\varkappa _{2}=x_{h}(T_{1},T_{3},e_{2}),$ $%
\varkappa _{3}=x_{h}(T_{2},T_{1},e_{3}),$ then, from Proposition \ref{dim
polyhedron}, the hyperbolic metric $h$ on $\partial \Delta $ is parametrized
by two of them, say $\varkappa _{1},$ $\varkappa _{2}.$ This implies that $%
\Delta $ is also parametrized by $\varkappa _{1},$ $\varkappa _{2}.$ It is
not difficult to express analytically $\alpha ,$ $\beta ,$ $\gamma $ as a
function of $\varkappa _{1},$ $\varkappa _{2},$ $\varkappa _{3}$ and
inversely. Therefore, we may derive the existence of a diffeomorphism $\phi :%
\mathbb{R}^{2}\rightarrow Int(T)$ \noindent which can be chosen to send $%
\varkappa _{1},$ $\varkappa _{2}$ to $\alpha ,$ $\beta .$ The expression of
the angles $\alpha ,$ $\beta ,$ $\gamma $ as a function of $\varkappa _{1},$ 
$\varkappa _{2},$ $\varkappa _{3}$ is indicated below.

Assume that $\Delta $ is projected to a Euclidean triangle $ABC$ in the $%
(x,y)-$plane and the angle at the vertex $A$ (resp. $B,$ $C)$ of $ABC$ is
equal to $\alpha $ (resp. $\beta ,$ $\gamma ),$ see Figure 1. Assuming,
without loss of generality, that the Euclidean length of $BC$ is equal to 1,
we have that

\begin{equation*}
|AB|=\frac{\sin \gamma }{\sin \alpha },\qquad |AC|=\frac{\sin \beta }{\sin
\alpha }
\end{equation*}

The shift parameters $\varkappa _{1},$ $\varkappa _{2},$ $\varkappa _{3}$
are given by the formulas:

\begin{equation*}
\varkappa _{1}=\log \frac{\sin \beta }{\sin \alpha }-\log \frac{\sin \gamma 
}{\sin \alpha },\text{ }\varkappa _{3}=\log 1-\log \frac{\sin \beta }{\sin
\alpha },\text{ }\varkappa _{2}=\log \frac{\sin \gamma }{\sin a}-\log 1
\end{equation*}

\noindent Therefore,

\begin{equation}
\varkappa _{1}=\log \frac{\sin \beta }{\sin \gamma },\text{ }\varkappa
_{2}=\log \frac{\sin \gamma }{\sin \alpha },\text{ }\varkappa _{3}=\log 
\frac{\sin \alpha }{\sin \beta }  \tag{3}  \label{three}
\end{equation}

\noindent Now we have

\begin{equation*}
e^{\varkappa _{2}}=\frac{\sin \gamma }{\sin \alpha }=\frac{\sin (\pi -\alpha
-\beta )}{\sin \alpha }=\cos \beta +\cos \alpha \frac{\sin \beta }{\sin
\alpha }\Longrightarrow
\end{equation*}

\begin{equation}
e^{\varkappa _{2}}=\cos \beta +\cos \alpha \cdot e^{-\varkappa _{3}}  \tag{4}
\label{four}
\end{equation}

\noindent From relation \ref{three} and replacing $\cos \beta $ from \ref
{four} we have that

\begin{equation*}
\cos \alpha =\frac{e^{2\varkappa _{2}}+e^{-2\varkappa _{3}}-1}{2e^{\varkappa
_{2}-\varkappa _{3}}}=\frac{e^{\varkappa _{2}+\varkappa _{3}}+e^{-\varkappa
_{2}-\varkappa _{3}}-e^{\varkappa _{3}-\varkappa _{2}}}{2}
\end{equation*}

\noindent Therefore

\begin{equation*}
\alpha =Arc\cos (\frac{e^{\varkappa _{2}+\varkappa _{3}}+e^{-\varkappa
_{2}-\varkappa _{3}}-e^{\varkappa _{3}-\varkappa _{2}}}{2}),
\end{equation*}

\noindent\ In a similar way we may express $\beta $ and $\gamma $ as a
function of $\varkappa _{2}$ and $\varkappa _{3}.$

Now, fix an ideal hyperbolic structure $h$ on $M$ and let $X=|\mathcal{D}%
^{(2)}|.$ Let $e_{i},$ $i=1,..,n$ be the edges of $M^{o}$ and let $\theta
_{h}(e_{i})$ be the angle around the edge $e_{i}.$ From the discussion above
the angles $\theta _{h}(e_{i})$ can be expressed as a function of shift
parameters of $X,$ but it is dificult to express the shift parameters as a
function of the angles $\theta _{h}(e_{i}).$ Therefore it is an interesting
problem to see, at least in some cases, whether or not $\theta _{h}(e_{i})$
determine the hyperbolic structure $h.$ In the next section we give examples
which explore this problem.

In the following proposition we investigate the linear relation among the
angles $\theta _{h}(e_{i}).$

\noindent

\begin{proposition}
\label{basic edges}Let \ $h$ be a hyperbolic metric on $M^{o}.$ Assuming
that $\dim \mathcal{T}_{\mathcal{D}}(M)=d,$ we may choose $d$ edges $e_{1},...,e_{d}$ such
that $\theta _{h}(e_{1}),...,\theta _{h}(e_{d})$ determine all $\theta
_{h}(e)$ for each edge $e.$
\end{proposition}

\begin{proof}
The space $M^{o}$ is obtained by gluing ideal hyperbolic tetrahedra by
isometries along their faces. Recall that $M^{o}$ does not have boundary,
i.e. each face of a tetrahedron is glued necessarily with another face.

Our lemma will be proved by induction on the number of pairs of faces which
are glued together in order to construct $M^{o}.$ For this reason, we are
obliged to prove the lemma in a more general context. We consider spaces $N$
which are constructed as follows:

(1) $N$ is always obtained by gluing ideal tetrahedra but we permit $N$ to
have free faces i.e. faces which is not glued to another face.

(2) $N$ is complete but not necessarily connected.

Let $v_{1},...,v_{l}$ be the cusps of $N$ and let $S_{h}^{1},...,S_{h}^{l}$
be the geometrical links of these cusps. Let $S_{h}=\cup _{i}S_{h}^{i}.$
Each edge $e$ of $\mathcal{D}$ intersects $S_{h}$ into two points which
appear as vertices of $S_{h}.$ Let $s$ be such an intersection point of $e$
with $S_{h}.$ If $\theta (s)$ denotes the angle around $s$ in $S_{h},$ then $%
\theta (s)=\theta _{h}(e).$ For each $i,$ $S_{h}^{i}$ is a compact Euclidean
surface with conical singularities, probably with boundary. Let $%
s_{i}^{1},...,s_{i}^{n_{i}},$ $i=1,2,..,l$ be the vertices of each $%
S_{h}^{i}.$ Then we may prove easily that

\begin{equation*}
\sum_{j=1,...,n_{i}}\phi _{i}^{j}=\pi \chi (DS_{h}^{i}),\forall
S_{h}^{i}\qquad (\ast )
\end{equation*}

\noindent \noindent \noindent where,

$\qquad \qquad \qquad \qquad \qquad \quad \phi _{i}^{j}=\left\{ 
\begin{array}{c}
\pi -\theta (s_{i}^{j}),\text{ if }s_{i}^{j}\in \partial S_{h}^{i} \\ 
2\pi -\theta (s_{i}^{j}),\text{ otherwise}
\end{array}
\right. $

\noindent \noindent and $DS_{h}^{i}$ is the double of the surface $S_{h}^{i}$
obtained by gluing two copies of $S_{h}^{i}$ along their boundaries provided
that $\partial S_{h}^{i}\neq \emptyset .$

Remark that for each variable $\phi _{i_{1}}^{j_{1}}$ in the set of
variables $\{\phi _{i}^{j}\}_{i,j}$ there exists exactly one variable $\phi
_{i_{2}}^{j_{2}}$ equal to $\phi _{i_{1}}^{j_{1}}.$ These two angles
correspond to vertices of $S_{h}=\cup _{i}S_{h}^{i}$ which are induced by
the same edge of $N.$

In order to prove our lemma it suffices to prove that the $l$ equations $%
(\ast )$ are linearly independent. Actually, the constants in the right hand
side of equations $(\ast )$ do not affect the linear independence of them.
So, it is sufficient to prove that the sums $\sum_{j=1,...,n_{i}}\phi
_{i}^{j},$ $i=1,..,l$ are linearly independent in the sense that we cannot
obtain one of them as a linear combination of the others. We will prove the
lemma by induction on the number of gluings between faces of ideal
tetrahedra.

Let $n=0.$ This means that we do not have any gluing between tetrahedra. So,
it suffices to examine the case where our space consists only of one
tetrahedron, because the variables that correspond to distinct tetrahedra
are distinct.

For an ideal hyperbolic tetrahedron $\Delta =v_{1}v_{2}v_{3}v_{4}$ the left
hand side of our equations are 
\begin{equation*}
\begin{array}{c}
\varphi _{1}^{1}+\varphi _{1}^{2}+\varphi _{1}^{3} \\ 
\varphi _{2}^{1}+\varphi _{2}^{2}+\varphi _{2}^{3} \\ 
\varphi _{3}^{1}+\varphi _{3}^{2}+\varphi _{3}^{3} \\ 
\varphi _{4}^{1}+\varphi _{4}^{2}+\varphi _{4}^{3}
\end{array}
\end{equation*}
Therefore, the matrix that correspond to these variables is 
\begin{equation*}
\left[ 
\begin{array}{cccccc}
1 & 1 & 1 & 0 & 0 & 0 \\ 
0 & 0 & 1 & 1 & 1 & 0 \\ 
0 & 1 & 0 & 1 & 0 & 1 \\ 
1 & 0 & 0 & 0 & 1 & 1
\end{array}
\right]
\end{equation*}
which has rank $4.$

Now, let $\{\sum_{j}\phi _{i}^{j}\}_{i}$ be the set of sums which correspond
to a space $N.$ For each $i$ we correspond a row and we assume that all rows
are linearly independent. We will show that if we glue two free faces of
tetrahedra so that the resulting space is complete the corresponding set of
sums, say $\{\sum_{m}\widetilde{\phi }_{n}^{m}\}_{n}$ is linearly
independent.

We remark that the set $\{\sum_{m}\widetilde{\phi }_{n}^{m}\}_{n}$ results
from the relations $\{\sum_{j}\phi _{i}^{j}\}_{i}$ by applying successively
the following two transformations:

\noindent$(A)$ Replace two rows by their sums. That is, 
\begin{equation*}
\{\sum_{j}\phi _{i}^{j}\}_{i}\rightarrow \{\sum_{j}\phi _{i}^{j}\}_{i\neq
i_{1},i_{2}}\cup \{\sum_{j}\phi _{i_{1}}^{j}+\sum_{j}\phi _{i_{2}}^{j}\}.
\end{equation*}

\noindent$(B)$ Replace two variables, by a new one, say $\phi _{0}.$ That
is, 
\begin{equation*}
\phi _{i_{1}}^{j_{1}}+\phi _{i_{1}}^{j_{2}}\rightarrow \phi _{0}
\end{equation*}

\noindent \noindent Remark that $\phi _{i_{1}}^{j_{1}}+\phi _{i_{1}}^{j_{2}}$
can appear either once in two different rows or twice in the same row.
Indeed, assuming that two free faces of $N$ are glued, then two different
edges, say $e^{\prime },$ $e^{\prime \prime },$ match together and give a
common edge $e.$ In the level of geometrical links of cusps we have the
following two possibilities:

\noindent $(i)$ two surfaces $S_{h}^{i_{1}},$ $S_{h}^{i_{2}}$ which are the
geometrical links of two distnct cusps of $N,$ are glued along an edge in
their boundaries or,

\noindent $(ii)$ two edges in the boundary of a surface $S_{h}^{i_{1}}$
which is the geometrical link of a cusp of $N,$ are glued together.

\noindent In both cases we obtain the relations $\{\sum_{m}\widetilde{\phi }%
_{n}^{m}\}_{n}$ from the relations $\{\sum_{j}\phi _{i}^{j}\}_{i}$ by
applying the rules $(A)$ and $(B).$

Obviously, every transformation of type $(A)$ or $(B)$ gives a system of
linearly independent sums and therefore our proof is complete.
\end{proof}

\noindent \includegraphics[scale=0.8]{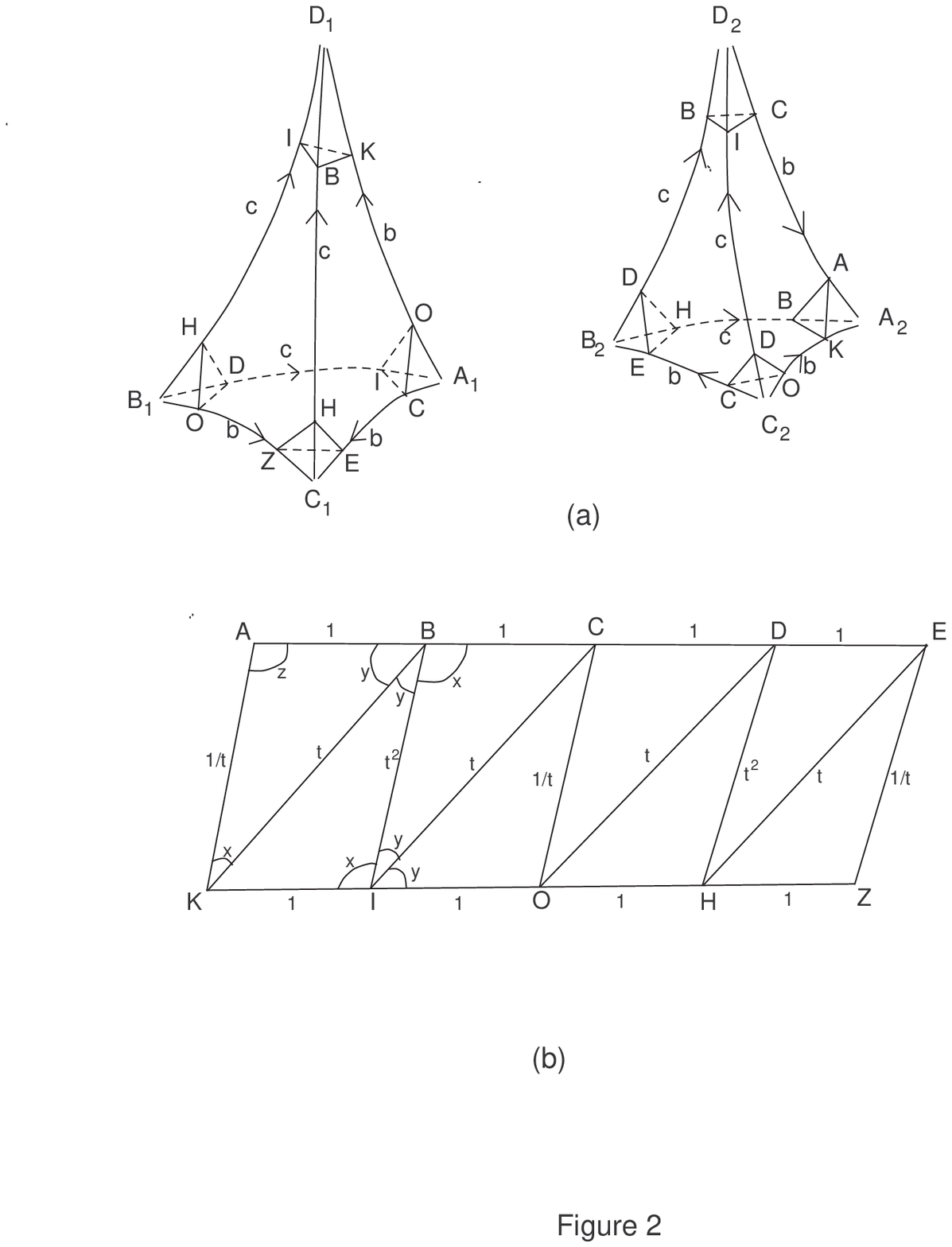}

\section{Examples}

Let $h\in \mathcal{T}_{\mathcal{D}}(M)$ and let $S_{h}$ be the geometrical link (not
necessarily connected) of $M^{o}$ with respect to $h.$ The fixed
triangulation $\mathcal{D}$ of $M$ induces on $S_{h}$ a (fixed) loose
triangulation\textit{\ }$\mathcal{D}_{e}.$ Every ideal hyperbolic
tetrahedron of $\mathcal{D}$ induces four similar pseudo-triangles on $%
S_{h}. $ Therefore $\mathcal{D}_{e}$ is equipped with a pattern which
indicates which Euclidean pseudo-triangles are pairwise similar. The
similarity of these pseudo-triangles imposes a system of equations between
the lengths of edges of $\mathcal{D}_{e}.$ The solutions of this system is
in $1-1$ correspondence with the elements $h\in \mathcal{T}_{\mathcal{D}}(M).$ This gives
rise to a parametrization of $\mathcal{T}_{\mathcal{D}}(M).$

The above idea is used to study $\mathcal{T}_{\mathcal{D}}(M)$ in the examples 2, 3 and 4
below.

\noindent \textbf{Example 1}

In the example 3.3.12 of \cite{[Thurstonbook]}, Thurston glues two ideal
hyperbolic tetrahedra so that the resulting space is a manifold $N$ with one
axis and one cusp. By truncating the tetrahedra we obtain a compact manifold 
$M$ with boundary and we have that $N$ is homeomorphic to $M^{o}.$ From
Theorem \ref{dimension}, the dimension of $\mathcal{T}_{\mathcal{D}}(M)$ is equal to $0.$
Therefore $M$ admits a unique ideal hyperbolic structure modulo the
equivalence relation of Definition \ref{defgemteichspace}. It is interesting
to remark here that if $\dim (\mathcal{T}_{\mathcal{D}}(M))=0,$ then the unique hyperbolic
structure in $\mathcal{T}_{\mathcal{D}}(M)$ is always obtained by gluing regular ideal
hyperbolic tetrahedra.

\noindent In \cite{[Fujii]}, all manifolds with one cusp and one axis are constructed
by gluing precisely two ideal hyperbolic tetrahedra. Using more than two
ideal tetrahedra it is easy to construct manifolds whose the number of edges
minus the number of cusps is zero. For all these manifolds, a similar
analysis and a similar result as above is valid.

\noindent \includegraphics[scale=0.8]{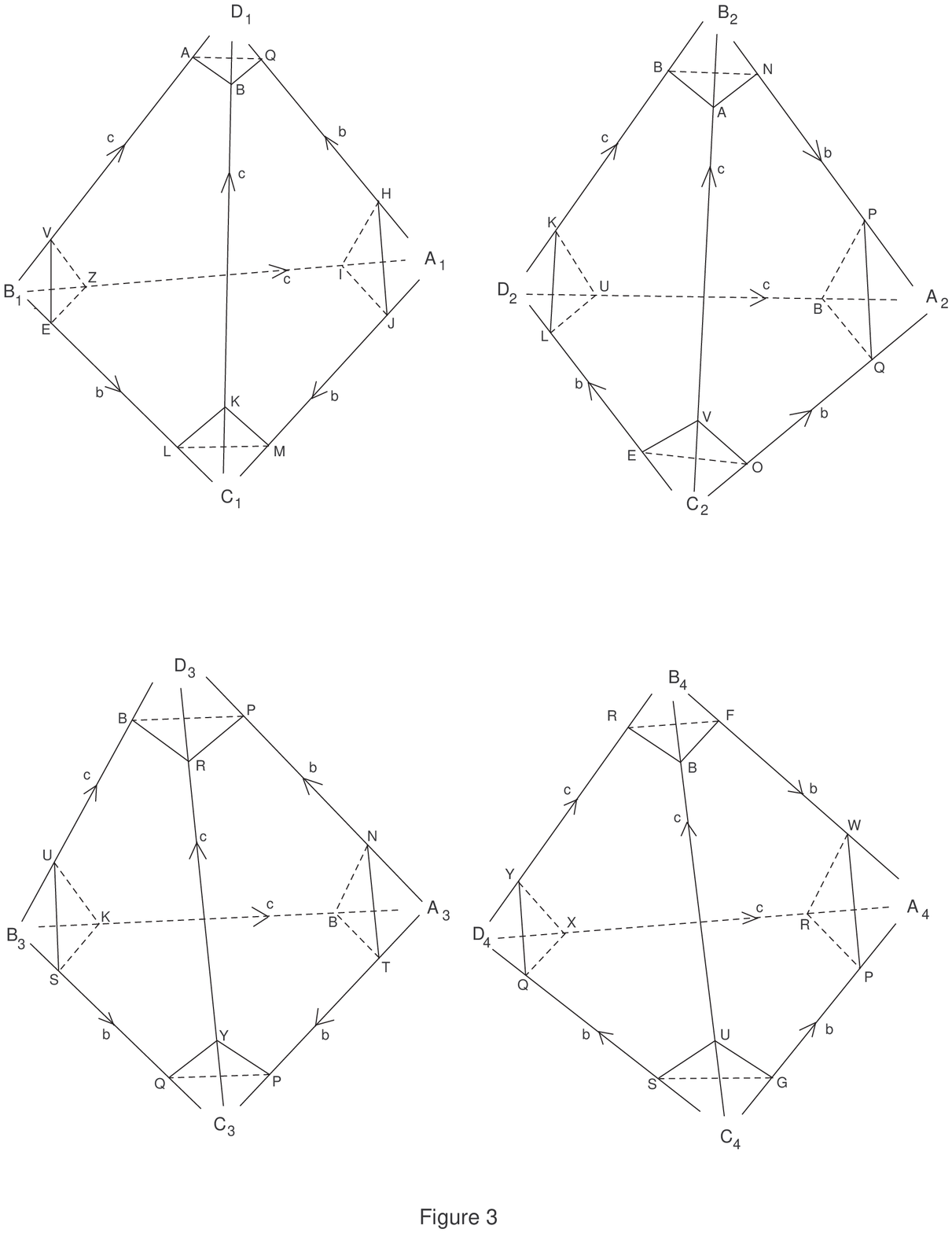}

\noindent \textbf{Example 2}

Consider two ideal hyperbolic tetrahedra $\Delta _{i}=A_{i}B_{i}C_{i}D_{i},$ 
$i=1,2,$ whose edges are labeled by the letters $b$ and $c$ and directed as
it is shown in Figure 2(a). There is a unique way to glue, via isometries,
the faces of $\Delta _{1},$ $\Delta _{2}$ so that the directed edges labeled
with $b$ (resp. $c)$ are identified to an edge, say $b$ (resp. $c).$ The
obtained manifold $N$ has one cusp and is homeomorphic to the interior $%
M^{o} $ of the compact manifold $M=\mathbb{S}^{3}-V,$ where $V$ is an open
tubular neighborhood of the figure eight knot in the sphere $\mathbb{S}^{3},$
see Example 1.4.8 in \cite{[Thurstonbook]}. From Theorem \ref{dimension},
the dimension of $T(M)$ is equal to $1.$

Fix an arbitrary ideal hyperbolic structure $h$ on $M$ and let $S_{h}$ be
the geometrical link of $v.$ The surface $S_{h}$ is obtained from the
polygon of Figure 2(b), by identifying the following directed segments:

\begin{center}
$AB\equiv KI,$ $BC\equiv IO,$ $CD\equiv OH,$ $DE\equiv HZ,$ $EZ\equiv AK.$ $%
\qquad \qquad $
\end{center}

We use the following notation: if $A,$ $B,$ $\Gamma $ are the vertices of a
Euclidean triangle $AB\Gamma $ and $|AB|,$ $|A\Gamma |,$ $|B\Gamma |$ the
lengths of its sides, then we set $|AB\Gamma |=$ $(|AB|,$ $|B\Gamma |,$ $%
|A\Gamma |).$ Remark now that all four triangles in $S_{h}$ induced by the
tetrahedron $\Delta _{1}$ (resp. $\Delta _{2})$ are similar. Therefore there
are positive numbers $\lambda ,$ $\mu ,$ $\nu ,$ $l,$ $m,$ $n$ such that we
have the following system of equations:

\begin{equation}
\begin{array}[c]{l}
|KBA| = \lambda |BIC|=\mu |ODC|=\nu |DHE|   \\
|IBK| = l|CIO|=m|HDO|=n|EHZ| 
\end{array}
\tag{5}  \label{five}
\end{equation}

We may assume that $|AB|=1$ and $|KB|=r.$ Then, using the equalities \ref
{five}, we may express successively all quantities in the system \ref{five}
as a function of the parameter $r,$ as follows:

$|AB|=|KI|=|BC|=|IO|=|CD|=|OH|=|DE|=|HZ|=m=\mu =1$

$|KB|=|IC|=|OD|=|HE|=l=n=\lambda =\nu =r$

$|BI|=|DH|=r^{2},$ $|AK|=|EZ|=|CO|=\frac{1}{r}$

Now, it is immediate to verify that the previous expressions of $r$ verify
all the equations of system \ref{five}. Also, from Theorem \ref{dimension},
we know that the dimension $d$ of $\mathcal{T}_{\mathcal{D}}(M)$ is equal to one. This
implies that $\mathcal{T}_{\mathcal{D}}(M)$ can be parametrized by the parameter $r.$

We also deduce that all triangles of $S_{h}$ are similar. We set $\widehat{%
AKB}=x,$ $\widehat{ABK}=y,$ $\widehat{BAK}=z$ see Figure 2(b). The
parameters $1,$ $r,$ $\frac{1}{r}$ are lengths of sides of a Euclidean
triangle. Therefore these quantities must satisfy the triangle inequalities
which imply that $\frac{\sqrt{5}-1}{2}<r<\frac{\sqrt{5}+1}{2}.$ If $r=1/r,$ $%
i.e.$ $r=1$ we deduce that all triangles of $S_{h}$ are equilateral and
therefore the tetrahedra $\Delta _{1}$ and $\Delta _{2}$ are regular and the
angles around the edges $b$ and $c$ are equal to $2\pi .$

Generally, the angle $\theta $ around the axis $c$ is equal to

\begin{center}
$\theta =2x+2y+2y$ $=2(\pi -z)+2y=$ $2\pi +2(y-z)$ \qquad \qquad $(\ast \ast
)$
\end{center}

We consider the function $\varphi (r)=z-y.$ Then, from the cosine law, we
have

\begin{equation*}
\varphi (r)=Arc\cos \frac{1+\frac{1}{r^{2}}-r^{2}}{2\frac{1}{r}}-Arc\cos 
\frac{1+r^{2}-\frac{1}{r^{2}}}{2r}=
\end{equation*}

\begin{equation*}
=Arc\cos \frac{1}{2}(r+\frac{1}{r}-r^{3})-Arc\cos \frac{1}{2}(\frac{1}{r}+r-%
\frac{1}{r^{3}})
\end{equation*}

We have that

\begin{eqnarray*}
\varphi ^{\prime }(r) &=&-\frac{1}{\sqrt{1-\frac{1}{4}(r+\frac{1}{r}%
-r^{3})^{2}}}\frac{1}{2}(1-\frac{1}{r^{2}}-3r^{2})+ \\
&&+\frac{1}{\sqrt{1-\frac{1}{4}(\frac{1}{r}+r-\frac{1}{r^{3}})^{2}}}\frac{1}{%
2}(-\frac{1}{r^{2}}+1+3r^{-4})
\end{eqnarray*}

But

\begin{equation*}
-1+\frac{1}{r^{2}}+3r^{2}=\frac{1}{r^{2}}(3r^{4}-r^{2}+1)>0
\end{equation*}

\begin{equation*}
-\frac{1}{r^{2}}+1+3r^{4}=\frac{1}{r^{4}}(r^{4}-r^{2}+3)>0
\end{equation*}

Therefore, $\varphi ^{\prime }(r)>0$ which implies that $\varphi $ is $1-1.$
Therefore, $r$ determines uniquely the angle $y-z$ and from equation $(\ast
\ast ),$ $r$ determines uniquely the angle $\theta $ around the edge $c.$
This proves that $\mathcal{T}_{\mathcal{D}}(M)$ is parametrized by the angle around the
edge $c$ of $N.\medskip $

\noindent \hspace*{-5cm} \includegraphics[scale=0.55]{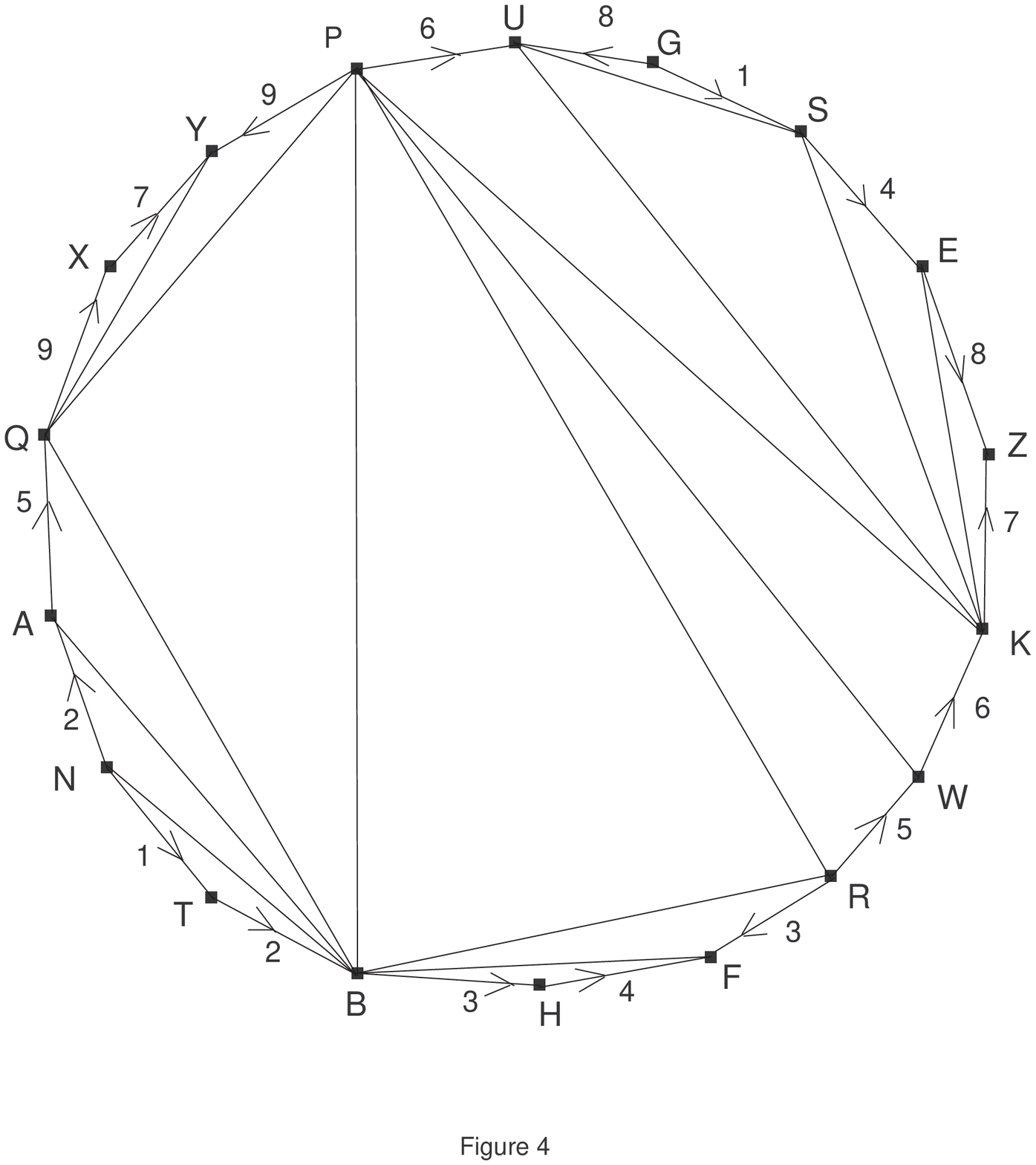}

\noindent \textbf{Example 3\medskip }

Consider four ideal hyperbolic tetrahedra $\Delta _{i}=A_{i}B_{i}C_{i}D_{i},$
$i=1,2,3,4$ whose edges are labeled by the letters $b$ and $c$ and directed
as it is shown in Figure 3. The faces of $\Delta _{i}$ are glued, via
isometries, so that the directed edges labeled with the same letter are
identified. More precisely the faces are glued as follows:\noindent

$B_{1}C_{1}D_{1}=C_{2}D_{2}B_{2},$ $A_{1}C_{1}D_{1}=C_{2}D_{2}A_{2},$ $%
D_{2}B_{2}A_{2}=B_{3}A_{3}D_{3},$

$C_{2}B_{2}A_{2}=B_{3}A_{3}C_{3},$ $B_{3}C_{3}D_{3}=C_{4}D_{4}B_{4},$ $%
A_{3}C_{3}D_{3}=C_{4}D_{4}A_{4},$

$D_{4}B_{4}A_{4}=$ $B_{1}A_{1}D_{1},$ $C_{4}B_{4}A_{4}=B_{1}A_{1}C_{1}.$

Denote by $N$ the obtained manifold and we remark that $N$ has two edges $b,$
$c$ and one cusp, say $v.$ By truncating the tetrahedra we obtain a compact
manifold $M$ with boundary. We have that $N$ is homeomorphic to $M^{o}.$ It
is not difficult to show that all the ideal structures on $M$ are hyperbolic
and we will prove that $\mathcal{T}_{\mathcal{D}}(M)$ is parametrized by the angle around
an edge of $N.$

Fix an arbitrary ideal hyperbolic structure $h$ on $M$ and let $S_{h}$ be
the geometrical link of $v.$ We may verify that $S_{h}$ is a closed surface
obtained from the polygon of Figure 4, by identifying its sides which are
labeled by the same number. By computing the Euler characteristic of $S_{h}$
we can see that $S_{h}$ is a surface of genus three.

$S_{h}$ contains four groups of four triangles which are similar because
they are induced by the same tetrahedron. Therefore, there are positive
numbers $\lambda ,\mu ,\nu ,$ $l,m,n,$ $r,s,t,$ $\rho ,\sigma ,\tau $ such
that:
\begin{equation}
\begin{array}[c]{l}
|IJH| =\lambda |KML|=\mu |ZVE|=\nu |BAQ|   \\
|ABN| =l|BQP|=m|UKL|=n|VOE|  \\
|BTN| =r|RBP|=s|YPQ|=t|KUS|   \\
|BRF| =\rho |RPW|=\sigma |UGS|=\tau |XYQ|  
\end{array}
\tag{6}  \label{six}
\end{equation}

\noindent On the other hand, due to the identifications of the faces of
tetrahedra we have also the following equalities:

\noindent $|MK|=|LU|,$ $|HJ|=|OE|,$ $|OV|=|SK|,$ $|NA|=|TB|,$ $|TN|=|SG|,$ $%
|PY|=|QX|,$ \noindent $|VZ|=|XY|,$ $|IH|=|RF|,$ $|AQ|=|RW|,$ $|LM|=|PW|,$ $%
|EZ|=|GU|,$ $|JI|=|FB|$

Without loss of generality, we may assume that $|IH|=1.$ Then, using the
above equalities, we may express successively all quantities in the system 
\ref{six} as a function of the parameter $r,$ as follows:

\noindent $|IH|=|KL|=|VE|=|AQ|=|BN|=|BP|=|YQ|=|US|=|RF|=|RW|=\tau =\lambda
=m=s=1$

\noindent $%
|IJ|=|KM|=|ZE|=|BQ|=|AN|=|UL|=|VO|=|BT|=|RP|=|YP|=|KS|=|BF|=|UG|=|XQ|=n=%
\sigma =l=\rho =1/r$

\noindent $|ZV|=|BA|=|UK|=|RB|=|XY|=1/r^{2}$

\noindent $|JH|=|ML|=|QP|=|OE|=|TN|=|PW|=|GS|=t=\nu =\mu =r$

Now, it is immediate to verify that the previous expressions of $r$ verify
all the equations of system \ref{six}. Also, from Theorem \ref{dimension},
we know that the dimension $d$ of $\mathcal{T}_{\mathcal{D}}(M)$ is equal to one. This
implies that $r$ parametrizes $\mathcal{T}_{\mathcal{D}}(M).$

The Eucliean triangles induced on $S_{h}$ are all similar triangles and the
lengths of their sides are either $(1,r,1/r)$ or $(1,1/r,1/r^{2}).$ Consider
the triangle $IJH,$ see Figure 3, and set $x=\widehat{IJH},$ $y=\widehat{JHI}%
,$ $z=\widehat{HIJ}.$ Therefore, all the triangles in $S_{h}$ have angles
equal to $x,$ $y,$ $z.$ Furthermore, we have that $%
|IJH|=|BQP|=|BTN|=|RPW|=(1/r,r,1).$ So, the angle aroud the axis $b$ is
easily computed and is equal to $4(x+2y)=4\pi +4(y-z).$

As is Example 2, we consider the function $\varphi (r)=z-y$ and we prove in
the same way that $\varphi ^{\prime }(r)>0$ and so $\varphi $ is $1-1.$
Therefore $\mathcal{T}_{\mathcal{D}}(M)$ is parametrized by the angle around the edge $%
b.\smallskip $

\includegraphics[scale=0.9]{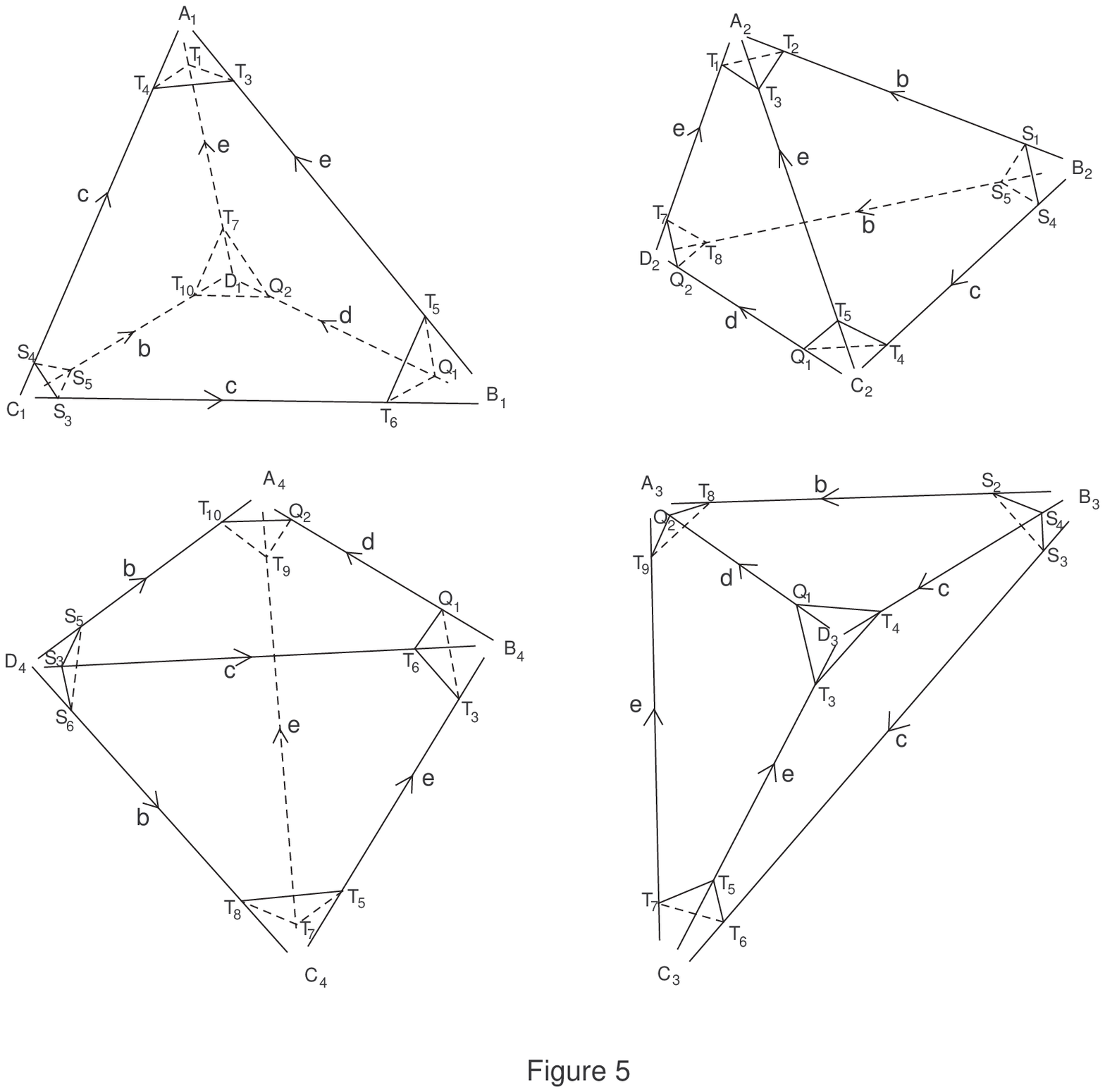}

\noindent \textbf{Example 4\medskip }

Consider four ideal hyperbolic tetrahedra $\Delta _{i}=A_{i}B_{i}C_{i}D_{i},$
$i=1,2,3,4$ whose edges are labeled by the letters $e,$ $b,$ $c,$ $d$ and
directed as it is shown in Figure 5. There is a unique way to glue, via
isometries, the faces of $\Delta _{i}$ so that the directed edges labeled
with the same letter are identified. The obtained manifold $N$ is
homeomorphic to the interior $M^{o}$ of the compact manifold $M=\mathbb{S}%
^{3}-V,$ where $V$ is an open tubular neighborhood of the Whitehead link in $%
\mathbb{S}^{3}$ (see p. 452 in \cite{[Ractliffe]}). Fix an arbitrary ideal
hyperbolic structure $h$ on $M.$ The manifold $N$ has two cusps, say $v,$ $%
w, $ and four edges so $\dim (\mathcal{T}_{\mathcal{D}}(M))=2.$ The geometrical link $%
S_{h}^{w}$ of $w$ (resp. $S_{h}^{v}$ of $v)$ is a Euclidean torus with
conical singularities, see Figure 6(a) (resp. Figure 6(b)). There are four
groups of four triangles which are similar because they are induced by the
same tetrahedron. Therefore, there are positive numbers $k_{i},$ $l_{i},$ $%
m_{i},$ $n_{i},$ $i=1,2,3$ such that:

\begin{equation}
\begin{array}[c]{l}
|T_{1}T_{3}T_{4}|=k_{1}|T_{7}T_{10}Q_{2}|=k_{2}|T_{6}T_{5}Q_{1}|=k_{3}|S_{3}S_{5}S_{4}| \\
|T_{1}T_{2}T_{3}|=l_{1}|T_{4}Q_{1}T_{5}|=l_{2}|S_{4}S_{1}S_{2}|=l_{3}|T_{7}Q_{2}T_{8}|  \\
|Q_{2}T_{9}T_{10}|=m_{1}|T_{8}T_{7}T_{5}|=m_{2}|Q_{1}T_{6}T_{3}|=m_{3}|S_{6}S_{3}S_{5}| \\
|Q_{2}T_{8}T_{9}|=n_{1}|Q_{1}T_{3}T_{4}|=n_{2}|S_{3}S_{2}S_{4}|=n_{3}|T_{6}T_{5}T_{7}| 
\end{array} \tag{7}  \label{seven}
\end{equation}

\noindent\includegraphics[scale=0.85]{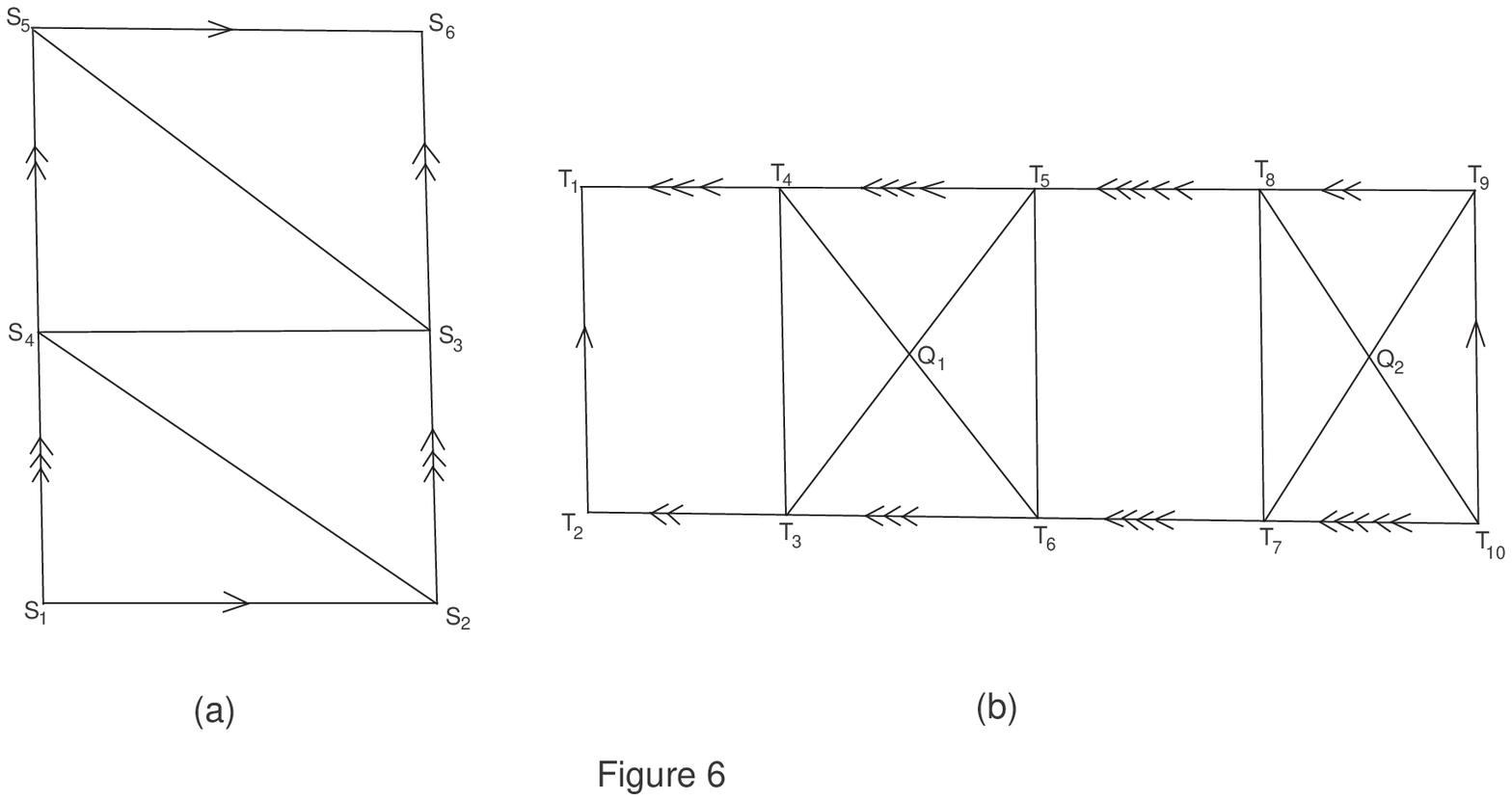}

Without loss of generality we assume that two edges in $S_{h}^{w}$ and $%
S_{h}^{v}$ respectively have length equal to $1.$ For example $%
|S_{4}S_{3}|=1,$ $|T_{3}T_{6}|=1.$ Using these two relations and the system
of equations \ref{seven} the lenghts of edges of $S_{h}=S_{h}^{w}\cup
S_{h}^{v}$ are computed as follows:

$%
|S_{4}S_{3}|=|T_{3}T_{6}|=|T_{1}T_{4}|=|T_{4}T_{5}|=|T_{6}T_{7}|=l_{2}=k_{3}=1 
$

$%
|T_{1}T_{2}|=|T_{9}T_{10}|=|T_{5}T_{6}|=|T_{8}T_{7}|=|S_{2}S_{3}|=|S_{1}S_{4}|=|T_{3}T_{4}|=|S_{3}S_{6}|= 
$

$|S_{4}S_{5}|=m_{2}=n_{1}=t$

$%
|T_{3}T_{2}|=|T_{9}T_{8}|=|T_{10}T_{7}|=|T_{8}T_{5}|=|S_{5}S_{6}|=|S_{1}S_{2}|=t^{2} 
$

$|S_{2}S_{4}|=|T_{1}T_{3}|=|T_{5}T_{7}|=|S_{5}S_{3}|=l_{1}=s$

$|T_{4}Q_{1}|=|Q_{1}T_{6}|=m_{3}=1/l_{3}=1/k_{2}=m_{1}=t/s$

$%
|T_{9}Q_{2}|=|Q_{1}T_{3}|=|Q_{2}T_{7}|=|T_{5}Q_{1}|=n_{3}=1/k_{1}=n_{2}=t^{2}/s 
$

$|T_{8}Q_{2}|=|Q_{2}T_{10}|=t^{3}/s$

Since the above expressions satisfy system \ref{seven} and $\dim (\mathcal{T}%
(M))=2$ we deduce that the parameters $t,$ $s$ parametrize $\mathcal{T}_{\mathcal{D}}(M).$

The angles $\theta _{h}(c),$ $\theta _{h}(d)$ around the axis $c$ and $d$
are respectively equal to

$\theta _{h}(c)=\widehat{T_{1}T_{4}T_{3}}+\widehat{T_{4}T_{1}T_{3}}+\widehat{%
T_{2}T_{1}T_{3}}+\widehat{T_{5}T_{7}T_{8}}+\widehat{T_{7}T_{6}T_{5}}+%
\widehat{T_{5}T_{7}T_{6}}=$

$2(\widehat{T_{1}T_{4}T_{3}}+\widehat{T_{4}T_{1}T_{3}}+\widehat{%
T_{2}T_{1}T_{3}})=2(\pi -\widehat{T_{1}T_{3}T_{4}}+\widehat{T_{2}T_{1}T_{3}}%
) $

$\theta _{h}(d)=\widehat{T_{1}T_{4}T_{3}}+\widehat{T_{1}T_{2}T_{3}}+\widehat{%
T_{5}T_{8}T_{7}}+\widehat{T_{5}T_{6}T_{7}}=2(\widehat{T_{1}T_{4}T_{3}}+%
\widehat{T_{1}T_{2}T_{3}})$

\noindent If we set $x=\widehat{T_{1}T_{4}T_{3}}+\widehat{T_{1}T_{2}T_{3}}$
and $y=\widehat{T_{2}T_{1}T_{3}}-\widehat{T_{1}T_{3}T_{4}}$ we have that $%
\theta _{h}(c)=2(\pi +y)$ and $\theta _{h}(d)=2x.$ Obviously $0<x<2\pi ,$ $%
-\pi <y<\pi $ and the angle around each edge of $N$ can be expressed as a
function of $x,$ $y.$ Now, in order to prove that $\mathcal{T}_{\mathcal{D}}(M)$ is
parametrized by the angles $x$ and $y$ it is sufficient to prove the
following claim.

\includegraphics[scale=1.13]{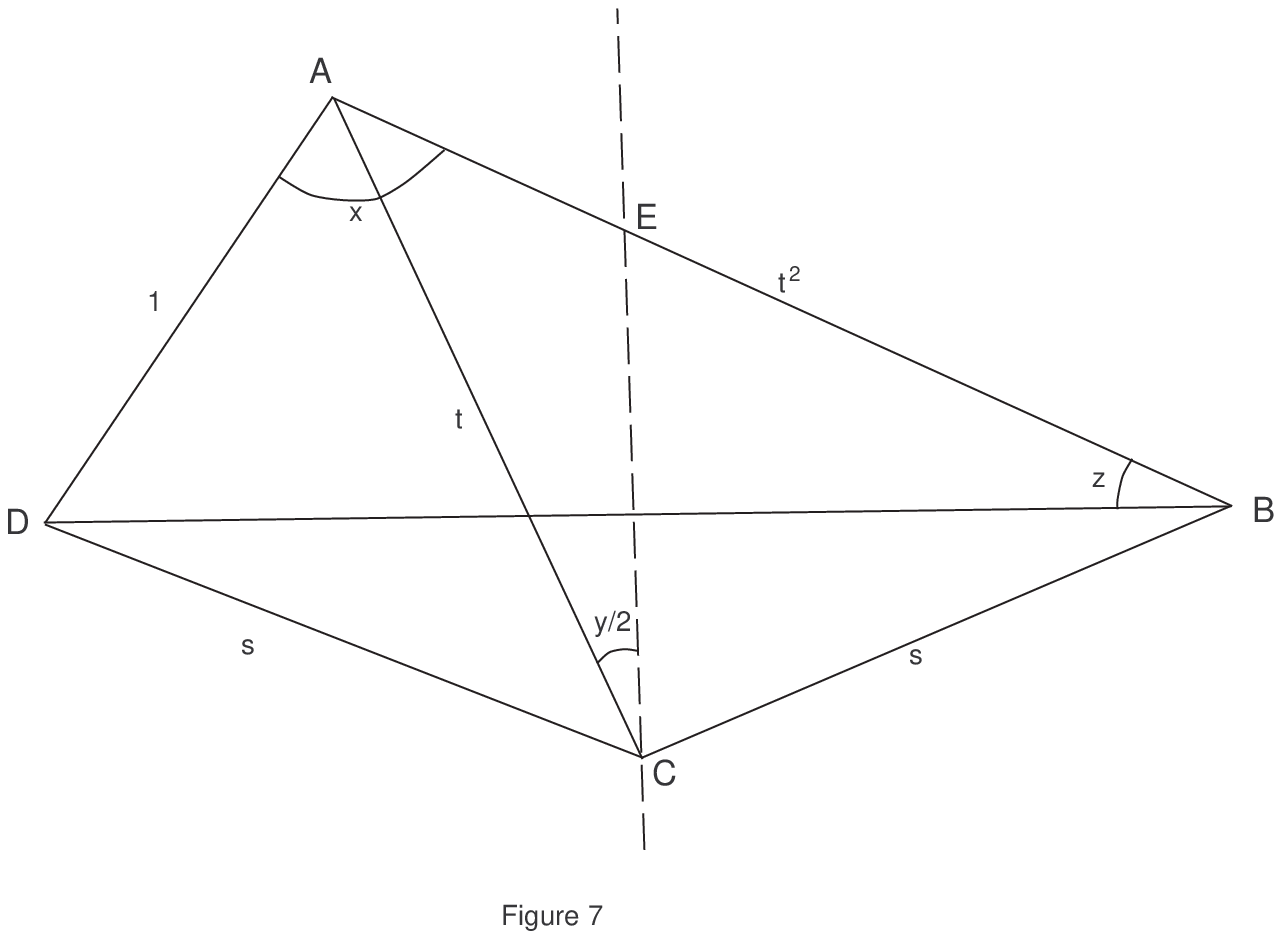}

\noindent \textit{Claim: }The mapping $\Theta :(t,$ $s)\rightarrow (x,$ $y)$
is $1-1.$

\noindent \textit{Proof of Claim}. From the expressions of edges of $S_{h}$
as a function of $t,$ $s$ we have that $|T_{1}T_{3}T_{4}|=(s,t,1)$ and $%
|T_{1}T_{2}T_{3}|=(t,t^{2},s).$ We glue the triangles $T_{1}T_{3}T_{4}$ and $%
T_{1}T_{2}T_{3}$ by identifying $T_{3}T_{4}$ with $T_{1}T_{2}.$ Thus, we
take a quadrilateral $ABCD$ such that $AC$ is a diagonal with $|AC|=t$ and $%
|DCA|=(s,t,1),$ $|CAB|=(t,t^{2},s),$ $x=\widehat{DAB},$ $y=\widehat{ACB}-%
\widehat{ACD},$ see Figure 7. Consider the perpendicular bisector $\zeta $
of $BD.$ Then $C\in \zeta $ and $AC$ and $\zeta $ form an angle equal to $%
\frac{y}{2}.$ Let $E=\zeta \cap AB.$

If $t\geq 1,$ we have that $y\geq 0.$ Let $z=\widehat{ABD},$ see Figure 7.
>From the \textit{law of sines }in the triangle $ACE$ we have that

\begin{equation*}
\frac{\sin \frac{y}{2}}{t^{2}-\frac{\sqrt{t^{4}+1-2t^{2}\cos x}}{2\cos z}}=%
\frac{\cos z}{t}
\end{equation*}

So,

\begin{equation*}
\sin \frac{y}{2}=\Phi (t,x)=t\sqrt{1-\frac{\sin ^{2}x}{t^{4}+1-2t^{2}\cos x}}%
-\frac{\sqrt{t^{4}+1-2t^{2}\cos x}}{2t}=
\end{equation*}

\begin{equation*}
=\frac{t(t^{2}-\cos x)}{\sqrt{t^{4}+1-2t^{2}\cos x}}-\frac{\sqrt{%
t^{4}+1-2t^{2}\cos x}}{2t}=\frac{1}{2}(t^{3}-\frac{1}{t})(t^{4}+1-2t^{2}\cos
x)^{-\frac{1}{2}}
\end{equation*}

\noindent Therefore,\smallskip

\noindent \noindent $\frac{\partial }{\partial t}(\sin \frac{y}{2})=\frac{1}{%
2}(3t^{2}+\frac{1}{t^{2}})(t^{4}+1-2t^{2}\cos x)^{-\frac{1}{2}}+$

\noindent $+(\frac{1}{t}-t^{3})(t^{4}+1-2t^{2}\cos x)^{-\frac{3}{2}}$%
\noindent $(t^{3}-t\cos x)=$

\noindent $=\frac{1}{2t^{2}}(t^{4}+1-2t^{2}\cos x)^{-\frac{3}{2}}$\noindent $%
\lbrack (3t^{4}+1)(t^{4}+1-2t^{2}\cos x)+2t^{2}(1-t^{4})(t^{2}-\cos x)]$

\noindent $=\frac{1}{2t^{2}}(t^{4}+1-2t^{2}\cos x)^{-\frac{3}{2}%
}(t^{8}-4t^{6}\cos x+6t^{4}-4t^{2}\cos x+1)\geq $

\noindent $\geq \frac{1}{2t^{2}}(t^{4}+1-2t^{2}\cos x)^{-\frac{3}{2}%
}(t^{8}-4t^{6}+6t^{4}-4t^{2}+1)=$

\noindent $=\frac{1}{2t^{2}}(t^{4}+1-2t^{2}\cos x)^{-\frac{3}{2}%
}(t^{2}-1)^{4}\geq 0.$

Therefore, we conclude that \noindent $\frac{\partial }{\partial t}(\sin 
\frac{y}{2})>0,$ $\forall t\geq 1.$ This implies that the mapping which
express $\sin \frac{y}{2}$ as function of $t$ is $1-1.$ On the other hand,
we have that $0\leq \frac{y}{2}\leq \frac{\pi }{2},$ therefore the mapping $%
f $ which express $y$ as a function of $t$ is $1-1.$

Now, we assume first that the angles $x,$ $y$ satisfy $0<x<2\pi $ and $0\leq
y<\pi .$ Then, using the fact that $f$ is $1-1,$ we may construct a unique
quadrilatel $ABCD,$ as well as, a point $E\in AB$ such that: $|AD|=1,$ $%
|CD|=|CB|,$ $\widehat{BAD}=x$ and $\widehat{ACE}=\frac{y}{2},$ see Figure 7.
Therefore the parameters $t=|AC|$ and $s=|CB|$ are uniquely determined from $%
x$ and $y.$ Therefore $\Theta $ is $1-1$ in this case.

If $t\leq 1$ then $-\pi <y\leq 0$ and we have that $\sin \frac{y}{2}=-\Phi (%
\frac{1}{t},x).$ An easy computation shows again that $\frac{\partial }{%
\partial t}(\sin \frac{y}{2})>0.$ So in this case, we also have that $t,$ $s$
are uniquely determined from $x$ and $y.$ Therefore we conclude that $\Theta 
$ is $1-1$ and the claim is proved.\medskip

Finally, the following question arises naturally:

\noindent \textbf{Question}. Let $M$ be an orientable compact manifold with
boundary equipped with an ideal triangulation $\mathcal{D}$ and assume that $%
\dim (\mathcal{T}_{\mathcal{D}}(M))\geq 1.$ Do there exist distinct elements $h,$ $%
h^{\prime }\in \mathcal{T}_{\mathcal{D}}(M)$ which have the same angles around the edges
of $\mathcal{D}?$

Note that that even if the ideal structures on $M$ are not hyperbolic, the
authors do not know examples where these structures are not uniquely
determined by the angles around the edges.\\[2mm]
\noindent {\bf Acknowledgement}
The authors would like to thank the referee for his remarks which improved
significantly the paper in many ways.


\begin{thebibliography}{99}
\bibitem{[ChAth]}  Ch. Charitos, A. Papadopoulos, \textit{The geometry of
ideal 2}$-$\textit{dimensional simplicial complexes,} Glasgow Math. J. 43,
39-66, 2001.

\bibitem{[ChAth1]}  Ch. Charitos, A. Papadopoulos, \textit{Hyperbolic
structures and measured foliations on 2}$-$\textit{dimensional complexes},
Monatsh. Math. 139, 1-17, 2003.

\bibitem{[fmp]}  R. Frigerio, B. Martelli, C. Petronio, \textit{Small
hyperbolic 3-manifolds with geodesic boundary, }Experiment Math. 13,
171-184, 2004.

\bibitem{[FP]}  R. Frigerio, C. Petronio, \textit{Constructions and
recognition of hyperbolic 3-manifolds with geodesic boundary}, Trans. Amer.
Math. Soc. 356, 3243-3282, 2004.

\bibitem{[Fujii]}  M. Fujii, Hyperbolic 3-Manifolds with Totally Geodesic
Boundary which are decomposed into Hyperbolic Truncated Tetrahedra, Tokyo J.
Math. 13(2), 355-373, 1990.

\bibitem{[Gromov]}  M. Gromov, \textit{Hyperbolic groups}, MSRI
Publications, 1988.

\bibitem{[CDP]}  M. Coornaert, T. Delzant, A. Papadopoulos, \textit{%
G\''{e}ometrie et th\''{e}orie des groupes,} Lecture Notes in Mathematics
1441, Springer-Verlag, 1990.

\bibitem{[Paulin]}  F. Paulin, \textit{Constructions of hyperbolic groups
via hyperbolization of polyhedra,} In group theory from a geometrical
viewpoint. ICTP, Trieste, Italy, 1990, E. Ghys and A. Haefliger eds, 1991.

\bibitem{[Ractliffe]}  J. Ratcliffe, \textit{Foundations of Hyperolic
Manifolds,} Springer, 1994.

\bibitem{[Thurstonbook]}  W. Thurston, \textit{Three-Dimensional Geometry
and Topology,} Edited by Silvio Levy, Princeton University Press, NJ 1997.

\bibitem{[ThurstonNotes]}  W. Thurston, \textit{The Geometry and Topology of
3}$-$\textit{Manifolds}, Lecture Notes, Princeton Univ. Princeton, NJ, 1979.

\bibitem{[Troyanov0]}  M. Troyanov, \textit{Les surfaces Euclidiennes \`{a}
singularit\''{e}s coniques,} L'' Enseignement Math. 32, 79-84, 1986.
\end{thebibliography}
\end{document}